\documentclass[11pt]{article}
\usepackage[a4paper,margin=1in]{geometry}
\usepackage{amsmath,amssymb,amsthm,mathtools}
\usepackage{microtype}
\usepackage{enumitem}
\usepackage[hidelinks,bookmarks=false]{hyperref}
\hypersetup{
  pdftitle={ell-abundant numbers and extremal divisor sums},
  pdfauthor={Ahmed Gaber},
  pdfsubject={Restricted divisor sums and abundance theory},
  pdfkeywords={restricted divisor sums, multiplicative functions, ell-full divisors, ell-abundant numbers, superabundant numbers, colossally abundant numbers}
}
\usepackage{booktabs}

\newtheorem{theorem}{Theorem}[section]
\newtheorem{lemma}[theorem]{Lemma}
\newtheorem{proposition}[theorem]{Proposition}
\newtheorem{corollary}[theorem]{Corollary}
\newtheorem{conjecture}[theorem]{Conjecture}
\theoremstyle{definition}
\newtheorem{definition}[theorem]{Definition}
\newtheorem{example}[theorem]{Example}
\theoremstyle{remark}
\newtheorem{remark}[theorem]{Remark}

\newcommand{\sigell}{\sigma_{(\ell)}}
\newcommand{\Aell}{A_{\ell}}
\newcommand{\calA}{\mathcal A}
\newcommand{\calS}{\mathcal S}
\newcommand{\calC}{\mathcal C}

\newcommand{\R}{\mathbb R}
\DeclareMathOperator{\rad}{rad}

\title{$\ell$-abundant numbers and extremal divisor sums}
\author{Ahmed Gaber\\
\small Mathematics Education Program, Faculty of Education and Arts, Sohar 3111, Oman\\
\small Department of Mathematics, Faculty of Science, Ain Shams University, Cairo, Egypt\\
\small \texttt{a.gaber@sci.asu.edu.eg}}
\date{}

\begin{document}
\maketitle

\begin{abstract}
For an integer $\ell\ge1$, we study the sum $\sigma_{(\ell)}(n)$ of the divisors of $n$ whose positive prime exponents are at least $\ell$.  The associated abundant, perfect, superabundant, and colossally abundant numbers recover the classical notions at $\ell=1$.  We prove strict nesting and infinitude of the abundant classes and derive obstructions to perfection.  Our results extend classical work of Ramanujan, Alaoglu, and Erd\H{o}s on superabundant prime-factor structure and colossally abundant exponent selection.  For $\ell\ge2$, we characterize all maximizers of $\sigma_{(\ell)}(n)/n^{1+\varepsilon}$ for every $\varepsilon>0$.  An explicit logarithmic formula determines the optimal positive exponents, while comparison with exponent $0$ resolves prime entry and all ties.  This characterization also yields infinitely many $\ell$-colossally abundant numbers.
\end{abstract}


\section{Introduction}

Classical abundance theory studies the arithmetic and extremal behavior of the divisor sum $\sigma(n)=\sum_{d\mid n}d$.  Its early development is documented in Dickson's work \cite{DicksonEven1913,DicksonFiniteness1913,DicksonHistory2005}, with subsequent advances concerning primitive abundant numbers and restrictions on prime factors \cite{AmatoEtAl2019,Avidon1996,Cohen1983,Cohen1984,Iannucci2005,Ivic1985}.  Davenport, Chowla, and Erd\H{o}s established foundational distribution results \cite{Chowla1934,Erdos1934,Davenport1933}; later work refined the density estimates and related counting problems \cite{Deleglise1998,JenningsEtAl2014,KlyveEtAl2019,KobayashiPollack2014,Kobayashi2016}.

The extremal theory developed by Ramanujan and Alaoglu--Erd\H{o}s \cite{AlaogluErdos1944,Ramanujan1915,Ramanujan1997} singles out superabundant and colossally abundant integers through unusually large divisor sums.  Its connection with the Riemann hypothesis is illustrated by Robin's equivalence between RH and $\sigma(n)<e^\gamma n\log\log n$ for every $n>5040$, where $\gamma$ is Euler's constant \cite{Robin1984}; related criteria and analogues are discussed in \cite{CaveneyEtAl2011,Nicolas2022,WashingtonYang2021}.

Fix an integer $\ell\ge1$, and let $D_\ell(n)$ be the set of $\ell$-full divisors of $n$, including $1$.  Thus every prime occurring in an admissible divisor has exponent at least $\ell$.  Writing $\rad(d)=\prod_{p\mid d}p$, with $\rad(1)=1$, we define
\[
\sigma_{(\ell)}(n)
=\sum_{\substack{d\mid n\\ \rad(d)^\ell\mid d}}d,
\qquad
A_\ell(n)=\frac{\sigma_{(\ell)}(n)}n.
\]
At $\ell=1$ these are the ordinary divisor sum and abundancy index.  For $\ell\ge2$, the restriction creates a local gap: exponents $1,\ldots,\ell-1$ contribute only the divisor $1$, whereas the divisor $p^\ell$ becomes admissible at exponent $\ell$.  Related work of Bansal, Jaindungarwal, and Maji studies generalized sums over perfect-power divisors \cite{BansalJaindungarwalMaji2023}.  Here we use the sum over all $\ell$-full divisors to develop an abundance theory.  The parameter changes the admissible divisors while retaining the threshold $2$, unlike classical $\alpha$-abundance \cite{Cohen1984}.

We prove that the sets $\mathcal A_\ell=\{n:A_\ell(n)>2\}$ form a strictly decreasing chain and are all infinite.  We exclude $\ell$-perfect prime powers and establish a parity obstruction to the equation $A_\ell(n)=2$.  For $\ell\ge2$, every nontrivial $\ell$-superabundant integer has consecutive-prime support and nonincreasing exponents, each at least $\ell$.  Theorem~\ref{thm:SAstructure} extends the prime-factor structure theorem of Alaoglu and Erd\H{o}s \cite[Theorem~1]{AlaogluErdos1944}, with the additional lower bound $\ell$ on every positive exponent.

For the $\ell$-colossally abundant problem, we maximize $F_{\ell,\varepsilon}(n)=\sigma_{(\ell)}(n)/n^{1+\varepsilon}$ for $\varepsilon>0$.  This extends Ramanujan's generalized super highly composite construction \cite{Ramanujan1997} and the exponent-selection formula of Alaoglu and Erd\H{o}s \cite[Theorem~10]{AlaogluErdos1944}.  Theorems~\ref{thm:positivebranch} and~\ref{thm:globalrule} determine the optimal positive exponents and their comparison with exponent $0$.  Our transition theory also extends to restricted divisor sums the classical description of critical parameters and tied maximizers presented by Caveney, Nicolas, and Sondow \cite{CaveneyEtAl2012}: adjacent positive-exponent transitions and prime-entry transitions must be treated separately.  The resulting characterization proves the infinitude of $\ell$-colossally abundant numbers.

Sections~2 and~3 develop the abundant and perfect theories.  Section~4 establishes the superabundant structure, Sections~5 and~6 give the colossally abundant characterization and transition theory, and Section~7 presents explicit examples.

\section{$\ell$-abundant numbers}

We begin with the local arithmetic of the restricted divisor sum.  The definitions are arranged so that the classical case is recovered without any change of normalization.

All divisors are positive, and the arguments of arithmetic functions range over positive integers.  We write $\nu_p(n)$ for the exponent of a prime $p$ in $n$, with $\nu_p(n)=0$ when $p\nmid n$.  The symbol $\log$ denotes the natural logarithm, and $\log_p x=(\log x)/(\log p)$.  Empty products have value $1$.

\begin{proposition}[Classical recovery and monotonicity in $\ell$]\label{prop:classicalrecovery}
For every positive integer $n$,
\[
\sigma_{(1)}(n)=\sigma(n).
\]
Moreover, if $1\le\ell_1<\ell_2$, then
\[
\sigma_{(\ell_2)}(n)\le \sigma_{(\ell_1)}(n).
\]
\end{proposition}

\begin{proof}
When $\ell=1$, every divisor $d\mid n$ is admissible, since every positive prime exponent is at least $1$.  Thus $D_1(n)$ is the full divisor set and $\sigma_{(1)}(n)=\sigma(n)$.  If $\ell_1<\ell_2$, every $\ell_2$-admissible divisor is automatically $\ell_1$-admissible, so
\[
D_{\ell_2}(n)\subseteq D_{\ell_1}(n).
\]
All summands are positive, so summing over the two divisor sets proves the asserted inequality.
\end{proof}

Dividing the inequality by $n$ gives
\[
A_{\ell_2}(n)\le A_{\ell_1}(n)
\qquad(1\le\ell_1<\ell_2).
\]
In particular, every $\ell$-abundant number is classically abundant.

\begin{definition}
Fix an integer $\ell\ge1$.  For a positive integer $n$, define
\[
D_\ell(n):=\{d\mid n:\rad(d)^\ell\mid d\},
\]
and
\[
\sigell(n):=\sum_{d\in D_\ell(n)}d.
\]
Equivalently, every prime occurring in an admissible divisor $d$ occurs there to exponent at least $\ell$.
\end{definition}

\begin{lemma}\label{lem:local}
For a prime $p$ and integer $a\ge0$,
\[
\sigell(p^a)=
\begin{cases}
1,&0\le a<\ell,\\[1mm]
1+p^\ell+p^{\ell+1}+\cdots+p^a,&a\ge\ell.
\end{cases}
\]
Moreover, $\sigell$ is multiplicative.
\end{lemma}

\begin{proof}
For $p^a$, the admissible exponents are $0$ and, when $a\ge\ell$, the integers $\ell,\ell+1,\ldots,a$.  This gives the local formula.

Now suppose $\gcd(m,n)=1$.  Every divisor $d$ of $mn$ has a unique representation $d=uv$ with $u\mid m$ and $v\mid n$.  Because the prime supports of $m$ and $n$ are disjoint,
\[
d\in D_\ell(mn)\quad\Longleftrightarrow\quad u\in D_\ell(m)\ \text{and}\ v\in D_\ell(n).
\]
Therefore
\[
\sigell(mn)=\sum_{u\in D_\ell(m)}\sum_{v\in D_\ell(n)}uv
=\sigell(m)\sigell(n).
\]
\end{proof}

Summing the geometric progression gives the equivalent formula
\[
\sigell(p^a)=\frac{p^{a+1}-p^\ell+p-1}{p-1}
\qquad(a\ge\ell).
\]

Thus, if $n=\prod_{i=1}^k p_i^{a_i}$, then
\[
\sigell(n)=\prod_{i=1}^k\sigell(p_i^{a_i}).
\]
Define
\[
\Aell(n):=\frac{\sigell(n)}{n}.
\]
Since both $\sigell$ and $n\mapsto n$ are multiplicative, $\Aell$ is multiplicative.

\begin{definition}
A positive integer $n$ is called $\ell$-abundant, $\ell$-perfect, or $\ell$-deficient according as
\[
\Aell(n)>2,\qquad \Aell(n)=2,\qquad \Aell(n)<2,
\]
respectively.  Let $\calA_\ell$ denote the set of $\ell$-abundant numbers.
\end{definition}

\begin{lemma}\label{lem:primepowermonotone}
For a prime $p$ and an integer $a\ge0$,
\[
\Aell(p^a)=
\begin{cases}
p^{-a},&0\le a<\ell,\\[1mm]
\displaystyle\frac{p}{p-1}-\frac{p^\ell-p+1}{(p-1)p^a},&a\ge\ell.
\end{cases}
\]
\end{lemma}

\begin{proof}
For $0\le a<\ell$, Lemma~\ref{lem:local} gives $\sigell(p^a)=1$, so $\Aell(p^a)=p^{-a}$.  For $a\ge\ell$, the geometric-sum formula gives
\[
\Aell(p^a)
=\frac{p^{a+1}-p^\ell+p-1}{(p-1)p^a}
=\frac{p}{p-1}-\frac{p^\ell-p+1}{(p-1)p^a}.
\]
\end{proof}

For fixed $p$, these formulas show that $\Aell(p^a)$ is strictly decreasing on the integers $0\le a<\ell$ and strictly increasing for $a\ge\ell$.  Indeed, the quantity $\frac{p^\ell-p+1}{(p-1)p^a}$ is positive and strictly decreases to $0$ as $a\to\infty$.  Consequently,
\[
\lim_{a\to\infty}\Aell(p^a)=\frac{p}{p-1}.
\]

\begin{proposition}\label{prop:primepower}
For every $\ell\ge1$, no prime power is $\ell$-abundant or $\ell$-perfect.
\end{proposition}

\begin{proof}
If $a<\ell$, then $\Aell(p^a)=p^{-a}\le1$.  If $a\ge\ell$, Lemma~\ref{lem:primepowermonotone} gives
\[
\Aell(p^a)<\frac{p}{p-1}\le2.
\]
For $p>2$ the upper bound is already $<2$, while for $p=2$ the limiting value $2$ is not attained at any finite exponent.  Hence $\Aell(p^a)<2$ in all cases.
\end{proof}

\begin{example}\label{ex:12}
The integer $12=2^2\cdot3$ is classically abundant because
\[
\frac{\sigma(12)}{12}=\frac{28}{12}>2.
\]
For every $\ell\ge2$, however, the exponent of $3$ is below $\ell$.  In particular, for $\ell=2$,
\[
\sigma_{(2)}(12)=\sigma_{(2)}(2^2)\sigma_{(2)}(3)=5\cdot1=5,
\qquad
A_2(12)=\frac5{12}<2.
\]
This example shows that the $\ell$-abundant classes differ from the classical abundant class for every $\ell\ge2$.
\end{example}

\begin{theorem}\label{thm:strictchain}
For every $\ell\ge1$,
\[
\calA_{\ell+1}\subsetneq\calA_\ell.
\]
\end{theorem}

\begin{proof}
Since every $(\ell+1)$-admissible divisor is $\ell$-admissible,
\[
D_{\ell+1}(n)\subseteq D_\ell(n),
\]
whence
\[
\sigma_{(\ell+1)}(n)\le\sigma_{(\ell)}(n)
\]
for every $n$.  Thus $\calA_{\ell+1}\subseteq\calA_\ell$.

To prove strictness, fix $\ell$ and consider
\[
n_k:=2^k3^k5^\ell,\qquad k\ge\ell+1.
\]
By Lemma~\ref{lem:primepowermonotone},
\[
\lim_{k\to\infty}A_\ell(n_k)
=2\cdot\frac32\cdot\left(1+5^{-\ell}\right)
=3\left(1+5^{-\ell}\right)>2.
\]
On the other hand, the exponent of $5$ is below $\ell+1$, so
\[
A_{\ell+1}(5^\ell)=5^{-\ell},
\]
and hence
\[
\lim_{k\to\infty}A_{\ell+1}(n_k)
=2\cdot\frac32\cdot5^{-\ell}
=\frac{3}{5^\ell}<2.
\]
Therefore, for all sufficiently large $k$, one has $A_\ell(n_k)>2$ and $A_{\ell+1}(n_k)<2$.  Hence $n_k\in\calA_\ell\setminus\calA_{\ell+1}$.
\end{proof}

Consequently,
\[
\calA_1\supsetneq\calA_2\supsetneq\calA_3\supsetneq\cdots.
\]

\begin{theorem}\label{thm:infabundant}
For every fixed $\ell\ge1$, there are infinitely many $\ell$-abundant numbers.
\end{theorem}

\begin{proof}
Let
\[
n_k=2^k3^k,\qquad k\ge\ell.
\]
For a prime $p$ and $k\ge\ell$, Lemma~\ref{lem:primepowermonotone} gives
\[
A_\ell(p^k)=\frac{p}{p-1}-\frac{p^\ell-p+1}{(p-1)p^k}.
\]
Consequently
\[
A_\ell(2^k)=2-\frac{2^\ell-1}{2^k},
\qquad
A_\ell(3^k)=\frac32-\frac{3^\ell-2}{2\cdot3^k},
\]
and therefore
\[
A_\ell(n_k)=
\left(2-\frac{2^\ell-1}{2^k}\right)
\left(\frac32-\frac{3^\ell-2}{2\cdot3^k}\right)
\longrightarrow3.
\]
Since the limit is strictly larger than $2$, there exists $k_0(\ell)$ such that $A_\ell(n_k)>2$ for every $k\ge k_0(\ell)$.  The integers $n_k$ are distinct, so there are infinitely many $\ell$-abundant numbers.
\end{proof}

\section{$\ell$-perfect numbers}

The classical case $\ell=1$ is the usual theory of perfect numbers.  In the restricted setting one asks whether the equation
\[
\sigma_{(\ell)}(n)=2n
\]
can hold when $\ell\ge2$.  The local jump at exponent $\ell$ makes this problem substantially different from the classical one.  We begin with two immediate observations.

\begin{remark}\label{rem:perfectbasic}
Let $\mathcal P_\ell:=\{n:A_\ell(n)=2\}$.  Then
\[
\mathcal P_\ell\cap\mathcal A_\ell=\varnothing.
\]
For every $\ell\ge2$,
\[
A_\ell(6)=\frac16.
\]
\end{remark}

\begin{proof}
The first assertion is immediate from the mutually exclusive conditions $A_\ell(n)=2$ and $A_\ell(n)>2$.  For $6=2\cdot3$ and $\ell\ge2$, both prime exponents are below $\ell$, so multiplicativity gives $\sigma_{(\ell)}(6)=1$.
\end{proof}

Since $6$ is classically perfect and $A_\ell(6)\ne2$ for $\ell\ge2$, we obtain
\[
\mathcal P_1\nsubseteq\mathcal P_\ell
\qquad(\ell\ge2).
\]

For $\ell\ge2$, Proposition~\ref{prop:primepower} already rules out all prime powers.  The following parity obstruction gives a further necessary condition.

\begin{theorem}[Parity obstruction]\label{thm:parity}
Let $\ell\ge2$ and write
\[
n=2^a\prod_{i=1}^k p_i^{b_i},
\]
where the $p_i$ are distinct odd primes, $a\ge0$, and $b_i\ge1$.  If $n$ is $\ell$-perfect, then for at least one index $i$,
\[
b_i\ge\ell\qquad\text{and}\qquad b_i\equiv\ell\pmod 2.
\]
\end{theorem}

\begin{proof}
If $n$ is $\ell$-perfect, then
\[
\sigell(n)=2n,
\]
so $\sigell(n)$ is even.  For the prime $2$, the local factor $\sigell(2^a)$ is always odd: it equals $1$ when $a<\ell$, and otherwise it is $1$ plus a sum of even integers.

For an odd prime $p$, if $b<\ell$, then $\sigell(p^b)=1$ is odd.  If $b\ge\ell$, then
\[
\sigell(p^b)=1+p^\ell+p^{\ell+1}+\cdots+p^b
\]
is a sum of exactly
\[
1+(b-\ell+1)=b-\ell+2
\]
odd terms.  This local factor is even exactly when $b-\ell+2$ is even, equivalently when $b\equiv\ell\pmod2$.

Because $\sigell(n)$ is the product of the local factors and is even, at least one odd-prime local factor must be even.  The stated condition follows.
\end{proof}

\begin{corollary}\label{cor:l2parity}
If $n=2^a m$, where $m$ is odd and every prime divisor of $m$ occurs to an odd exponent, then $n$ is not $2$-perfect.
\end{corollary}

\begin{proof}
For $\ell=2$, Theorem~\ref{thm:parity} requires at least one odd-prime exponent to be even.  This contradicts the hypothesis.
\end{proof}

\begin{remark}[Odd $2$-abundant integers]\label{rem:oddabundant}
Odd $2$-abundant numbers occur.  For example,
\[
n=3^3 5^3 7^3 11^3
\]
satisfies
\[
A_2(n)
=\frac{37}{27}\frac{151}{125}\frac{393}{343}\frac{1453}{1331}
\approx2.0705746>2.
\]
Thus a universal inequality $A_2(n)<2$ for odd $n$ cannot be used to exclude odd $2$-perfect numbers.  The parity obstruction above provides only a necessary condition.
\end{remark}

\begin{conjecture}\label{conj:noperfect}
For every $\ell\ge2$, there are no $\ell$-perfect numbers.
\end{conjecture}

\section{$\ell$-superabundant numbers}

Classically, a positive integer $n$ is superabundant if $\sigma(m)/m<\sigma(n)/n$ for every positive integer $m<n$.  Such numbers have highly rigid prime factorizations and play an important role in the extremal theory of the divisor sum; see Alaoglu--Erd\H{o}s \cite{AlaogluErdos1944}, Erd\H{o}s--Nicolas \cite{ErdosNicolas1975}, Akbary--Friggstad \cite{AkbaryFriggstad2009}, and Caveney--Nicolas--Sondow \cite{CaveneyEtAl2012}.  We now show that the same two classical structural principles---consecutive prime support and nonincreasing exponents---survive for $A_\ell$, together with a new $\ell$-fullness condition.

\begin{definition}
A positive integer $n$ is called $\ell$-superabundant ($\ell$-SA) if
\[
A_\ell(m)<A_\ell(n)\qquad\text{for every }m<n.
\]
Let $\calS_\ell$ denote the set of $\ell$-superabundant numbers.
\end{definition}

\begin{definition}
A positive integer $n$ is called $\ell$-full if every prime divisor $p\mid n$ satisfies $\nu_p(n)\ge\ell$.
\end{definition}

\begin{theorem}\label{thm:classicalnoncomp}
For every $\ell\ge2$, neither $\calS_\ell\subseteq\calS_1$ nor $\calS_1\subseteq\calS_\ell$ holds.
\end{theorem}

\begin{proof}
The integer $2$ is classically superabundant, since $A_1(2)=3/2>A_1(1)=1$.  For every $\ell\ge2$,
\[
A_\ell(2)=\frac12<A_\ell(1)=1,
\]
so $2\notin\calS_\ell$.  Hence $\calS_1\nsubseteq\calS_\ell$.

For $\ell=2$,
\[
(A_2(1),\ldots,A_2(8))=\left(1,\frac12,\frac13,\frac54,\frac15,\frac16,\frac17,\frac{13}{8}\right),
\]
so $8\in\calS_2$.  Classically $A_1(6)=2>A_1(8)=15/8$, hence $8\notin\calS_1$.

Now let $\ell\ge3$.  If $m<2^\ell$, then no prime can occur in $m$ to exponent $\ell$ or more, because $2^\ell$ is the smallest such prime power.  Hence $\sigma_{(\ell)}(m)=1$ and $A_\ell(m)=1/m\le1$.  But
\[
A_\ell(2^\ell)=1+2^{-\ell}>1,
\]
so $2^\ell\in\calS_\ell$.  On the other hand, $6<2^\ell$ and
\[
A_1(6)=2>2-2^{-\ell}=A_1(2^\ell),
\]
so $2^\ell\notin\calS_1$.
\end{proof}

\begin{theorem}[No monotonicity in the parameter]\label{thm:S2S3}
Neither $\calS_2\subseteq\calS_3$ nor $\calS_3\subseteq\calS_2$.
\end{theorem}

\begin{proof}
First, $4\in\calS_2$ because $A_2(4)=5/4$ exceeds the values at $1,2,3$.  But
\[
A_3(4)=\frac14<A_3(1)=1,
\]
so $4\notin\calS_3$.

For the reverse inclusion, we prove that $256=2^8$ is $3$-superabundant.  Its value is
\[
A_3(256)=\frac{1+8+16+32+64+128+256}{256}=\frac{505}{256}.
\]
Let $m<256$.  If no prime exponent in $m$ is at least $3$, then $A_3(m)=1/m\le1$.  If exactly one prime, say $p$, occurs to exponent at least $3$, write $m=p^a u$, where every prime exponent in $u$ is below $3$.  Then
\[
A_3(m)=\frac{A_3(p^a)}{u}\le A_3(p^a).
\]
If $p=2$, then $a\le7$, so $A_3(p^a)\le A_3(2^7)=249/128$.  If $p\ge3$, then $a\le5$, and Lemmas~\ref{lem:primepowermonotone} and \ref{lem:pmonotone} below give
\[
A_3(p^a)\le A_3(3^5)=\frac{352}{243}<\frac{249}{128}.
\]
If at least two primes occur to exponent at least $3$, then $m\ge2^3 3^3=216$.  Any pair other than $2,3$ contributes at least $2^3 5^3=1000$.  Thus $216\mid m$, and its only positive multiple below $256$ is $m=216$, for which
\[
A_3(216)=\frac98\frac{28}{27}=\frac76.
\]
Thus
\[
A_3(m)\le\frac{249}{128}=\frac{498}{256}<\frac{505}{256}=A_3(256)
\]
for every $m<256$, proving $256\in\calS_3$.

Finally,
\[
A_2(256)=\frac{509}{256},
\qquad
A_2(144)=A_2(2^4)A_2(3^2)=\frac{29}{16}\frac{10}{9}=\frac{290}{144}.
\]
The inequality $290/144>509/256$, together with $144<256$, shows that $256\notin\calS_2$.
\end{proof}

\begin{lemma}\label{lem:pmonotone}
Fix $\ell\ge2$.
\begin{enumerate}[label=\textup{(\roman*)}]
\item For every integer $r\ge\ell$, the function $p\mapsto A_\ell(p^r)$ is strictly decreasing for real $p>1$, using the prime-power formula to define it for real $p$.
\item If $r>s\ge\ell$ are integers, then, as $p$ ranges over the primes,
\[
\frac{A_\ell(p^r)}{A_\ell(p^s)}
\]
is strictly decreasing with $p$.
\end{enumerate}
\end{lemma}

\begin{proof}
For (i), write $x=1/p$.  From Lemma~\ref{lem:local},
\[
A_\ell(p^r)=\sum_{j=0}^{r-\ell}x^j+x^r.
\]
This is strictly increasing in $x\in(0,1)$, hence strictly decreasing in $p$.

For (ii), it is enough to prove that for each $a\ge\ell$,
\[
q_a(p):=\frac{A_\ell(p^{a+1})}{A_\ell(p^a)}
\]
is strictly decreasing, because
\[
\frac{A_\ell(p^r)}{A_\ell(p^s)}=\prod_{a=s}^{r-1}q_a(p).
\]
Put
\[
C=C_\ell(p):=p^\ell-p+1>0.
\]
Using
\[
\sigma_{(\ell)}(p^a)=\frac{p^{a+1}-C}{p-1},
\]
we obtain
\[
q_a(p)=\frac1p\frac{p^{a+2}-C}{p^{a+1}-C}
=1+D_a(p),
\]
where
\[
D_a(p)=\frac{C(p-1)}{p\bigl(p^{a+1}-C\bigr)}>0.
\]
We first prove $D_a'(p)<0$ for real $p\ge3$, and then compare $p=2$ with $p=3$.  Since $D_a>0$, $D_a'$ has the sign of the logarithmic derivative.  A direct differentiation gives
\[
\frac{D_a'}{D_a}
=\frac{C'}C+\frac1{p-1}-\frac1p
-\frac{(a+1)p^a-C'}{p^{a+1}-C}.
\]
Multiplication by the positive factor $(p^{a+1}-C)/p^a$ shows that this has the same sign as
\[
E(p):=\frac{pC'}C-(a+1)+\frac{1-C/p^{a+1}}{p-1}.
\]
Because $C'=\ell p^{\ell-1}-1$, for $p\ge3$ we have
\[
E(p)
<\frac{\ell p^\ell-p}{C}-(\ell+1)+\frac1{p-1}
=-\frac{N_{\ell}(p)}{(p-1)C},
\]
where
\[
N_\ell(p)=(p-2)p^\ell-\ell(p-1)^2+2(p-1).
\]
For $\ell=2$,
\[
N_2(p)=p^3-4p^2+6p-4>0\qquad(p\ge3),
\]
because this polynomial is positive at $p=3$ and has positive derivative for $p\ge3$.  Moreover,
\[
N_{\ell+1}(p)-N_\ell(p)
=(p-1)\bigl((p-2)p^\ell-(p-1)\bigr)>0
\]
for $p\ge3$ and $\ell\ge2$.  Hence $E(p)<0$ for every $p\ge3$.

Thus $D_a(p)$ is strictly decreasing for real $p\ge3$.  To include the prime $2$, it remains only to prove
\[
D_a(2)>D_a(3).
\]
Set
\[
\alpha=2^\ell-1,\qquad \beta=3^\ell-2.
\]
Then
\[
\frac{D_a(2)}{D_a(3)}
=\frac{3\alpha}{4\beta}
\frac{3^{a+1}-\beta}{2^{a+1}-\alpha}.
\]
The last fraction is strictly increasing in $a\ge\ell$, because
\begin{align*}
&(3^{a+2}-\beta)(2^{a+1}-\alpha)
 -(3^{a+1}-\beta)(2^{a+2}-\alpha)\\
&\qquad =6^{a+1}-2\alpha3^{a+1}+\beta2^{a+1}\\
&\qquad =3^{a+1}(2^{a+1}-2\alpha)+\beta2^{a+1}>0,
\end{align*}
since $2^{a+1}-2\alpha\ge2$.  Hence it suffices to check $a=\ell$.  There
\[
D_\ell(2)=\frac{2^\ell-1}{2(2^\ell+1)},
\qquad
D_\ell(3)=\frac{3^\ell-2}{3(3^\ell+1)}.
\]
The inequality $D_\ell(2)>D_\ell(3)$ is equivalent to
\[
6^\ell+7\cdot2^\ell-5\cdot3^\ell+1>0.
\]
For $\ell=2$ the left side is $20$.  For $\ell\ge3$, already $6^\ell>5\cdot3^\ell$, so the inequality is immediate.  Therefore $q_a(2)>q_a(3)$, and together with the strict decrease on $[3,\infty)$ this proves (ii) on the primes.
\end{proof}

In particular, part~\textup{(i)} holds when $p$ is restricted to prime values.

\begin{theorem}[Structure of $\ell$-superabundant numbers]\label{thm:SAstructure}
Let $\ell\ge2$, let $n>1$ be $\ell$-superabundant, and write
\[
n=\prod_{i=1}^k p_i^{a_i},
\qquad p_1<p_2<\cdots<p_k,
\]
where the $p_i$ are distinct primes and the $a_i$ are positive integers.  Then:
\begin{enumerate}[label=\textup{(\roman*)}]
\item $a_i\ge\ell$ for every $i$;
\item every prime $q\le p_k$ divides $n$;
\item $a_1\ge a_2\ge\cdots\ge a_k$.
\end{enumerate}
\end{theorem}

\begin{proof}
For (i), suppose $p^a\Vert n$ with $1\le a<\ell$, and write $n=p^a m$, $(p,m)=1$.  Since $\sigma_{(\ell)}(p^a)=1$,
\[
A_\ell(n)=p^{-a}A_\ell(m)<A_\ell(m),
\]
while $m<n$, contradicting superabundance.

For (ii), suppose that a prime $q<p_k$ does not divide $n$.  Write $n=p_k^{a_k}N$ and set
\[
m=q^{a_k}N.
\]
Then $m<n$.  By (i), $a_k\ge\ell$, and Lemma~\ref{lem:pmonotone}(i) gives
\[
A_\ell(q^{a_k})>A_\ell(p_k^{a_k}).
\]
By multiplicativity, $A_\ell(m)>A_\ell(n)$, a contradiction.

For (iii), suppose $p_i<p_j$ but $a_i<a_j$.  By (i), both exponents are at least $\ell$.  Replace the two prime powers by swapping their exponents:
\[
m=p_i^{a_j}p_j^{a_i}\frac{n}{p_i^{a_i}p_j^{a_j}}.
\]
Since $p_i<p_j$ and $a_j>a_i$, we have $m<n$.  Moreover,
\[
\frac{A_\ell(m)}{A_\ell(n)}
=
\frac{A_\ell(p_i^{a_j})/A_\ell(p_i^{a_i})}
     {A_\ell(p_j^{a_j})/A_\ell(p_j^{a_i})}>1
\]
by Lemma~\ref{lem:pmonotone}(ii).  This again contradicts superabundance.  Hence the exponents are nonincreasing.
\end{proof}

In particular, $n$ is $\ell$-full and $p_1=2$.

\begin{remark}[Comparison with the classical structure]\label{rem:SAclassical}
The classical consecutive-prime support and exponent ordering are due to Alaoglu and Erd\H{o}s \cite[Theorem~1]{AlaogluErdos1944}.  Theorem~\ref{thm:SAstructure} extends this structure to $\ell\ge2$ and adds the restriction that every positive exponent is at least $\ell$.
\end{remark}

\begin{lemma}\label{lem:successiveSA}
Let $\ell\ge2$, and let $1<n_1<n_2$ be successive $\ell$-superabundant numbers.  Then
\[
n_2\le2n_1.
\]
\end{lemma}

\begin{proof}
By Theorem~\ref{thm:SAstructure}, write $n_1=2^aM$ with $a\ge\ell$ and $M$ odd.  Lemma~\ref{lem:primepowermonotone} gives
\[
A_\ell(2^{a+1})>A_\ell(2^a),
\]
so by multiplicativity
\[
A_\ell(2n_1)>A_\ell(n_1).
\]
Let $m_0$ be the least integer greater than $n_1$ such that $A_\ell(m_0)>A_\ell(n_1)$.  Then $m_0\le2n_1$, and the superabundance of $n_1$ together with the minimality of $m_0$ gives $A_\ell(m)<A_\ell(m_0)$ for every $m<m_0$.  Thus $m_0$ is $\ell$-superabundant, and $n_2\le m_0\le2n_1$.
\end{proof}

\section{$\ell$-colossally abundant numbers}

Ramanujan's generalized highly composite framework, including the material preserved in the continuation of his work and later edited by Nicolas and Robin, and the subsequent colossally abundant theory of Alaoglu and Erd\H{o}s are based on extremal quotients of the form $\sigma(n)/n^{1+\varepsilon}$ \cite{AlaogluErdos1944,Ramanujan1915,Ramanujan1997}.  Their importance in the maximal-order theory of $\sigma$ and in questions surrounding Robin's inequality has led to extensive subsequent work \cite{Briggs2006,CaveneyEtAl2011,CaveneyEtAl2012,Musin2020,NazardonyaviYakubovich2012,NazardonyaviYakubovich2014,Nicolas2022}.  We now formulate the exact restricted analogue.  The point of this section is not only existence: because $\sigma_{(\ell)}$ is multiplicative, the global extremal problem can be solved completely in terms of independent local prime-power optimizations.

\begin{definition}
Let $\ell\ge1$.  A positive integer $N$ is called $\ell$-colossally abundant ($\ell$-CA) if there exists $\varepsilon>0$ such that
\[
\frac{\sigma_{(\ell)}(N)}{N^{1+\varepsilon}}
\ge
\frac{\sigma_{(\ell)}(m)}{m^{1+\varepsilon}}
\qquad\text{for every }m\ge1.
\]
Let $\calC_\ell$ denote the set of all such integers.
\end{definition}

\begin{remark}[The classical case]
For $\ell=1$, the definition becomes the usual colossally abundant maximization problem with all tied maximizers retained; see Caveney, Nicolas, and Sondow \cite{CaveneyEtAl2012}.  The exponent-selection rule is compared explicitly with that of Alaoglu and Erd\H{o}s after Proposition~\ref{prop:l1T}.
\end{remark}

For fixed $\ell$ and $\varepsilon>0$, put
\[
F_{\ell,\varepsilon}(n):=\frac{\sigma_{(\ell)}(n)}{n^{1+\varepsilon}}.
\]

\begin{lemma}\label{lem:CAsubsetSA}
For every $\ell\ge1$,
\[
\calC_\ell\subseteq\calS_\ell.
\]
\end{lemma}

\begin{proof}
Let $N$ maximize $F_{\ell,\varepsilon}$ for some $\varepsilon>0$.  If $m<N$, then
\[
A_\ell(m)
=F_{\ell,\varepsilon}(m)m^\varepsilon
\le F_{\ell,\varepsilon}(N)m^\varepsilon
<F_{\ell,\varepsilon}(N)N^\varepsilon
=A_\ell(N).
\]
Hence $N$ is $\ell$-superabundant.
\end{proof}

The case $\ell=1$ is the classical inclusion of colossally abundant numbers among superabundant numbers \cite{AlaogluErdos1944}.

For a prime $p$ and integer $a\ge0$, define
\[
f_{p,\varepsilon}(a)
:=\frac{\sigma_{(\ell)}(p^a)}{p^{a(1+\varepsilon)}}.
\]
Then $f_{p,\varepsilon}(0)=1$, and multiplicativity gives
\[
F_{\ell,\varepsilon}(n)
=\prod_p f_{p,\varepsilon}(\nu_p(n)).
\]
For $a\ge\ell-1$, define
\[
R_\ell(p,a):=
\frac{\sigma_{(\ell)}(p^{a+1})}{\sigma_{(\ell)}(p^a)},
\qquad
G_\ell(p,a):=\log_p R_\ell(p,a).
\]

\begin{lemma}\label{lem:Rdecreasing}
Fix $\ell\ge2$ and a prime $p$.  The sequence $R_\ell(p,a)$ is strictly decreasing for integers $a\ge\ell-1$.  Moreover,
\[
R_\ell(p,a)>p\qquad(a\ge\ell-1),
\]
and
\[
\lim_{a\to\infty}R_\ell(p,a)=p.
\]
\end{lemma}

\begin{proof}
For $a=\ell-1$,
\[
R_\ell(p,\ell-1)=1+p^\ell>p.
\]
For $a\ge\ell$,
\[
R_\ell(p,a)
=1+\frac{p^{a+1}}{\sigma_{(\ell)}(p^a)}.
\]
The inequality $R_\ell(p,a)>R_\ell(p,a+1)$ is equivalent to
\[
p^{a+1}>(p-1)\sigma_{(\ell)}(p^a).
\]
Using Lemma~\ref{lem:local}, the right side equals
\[
p^{a+1}-p^\ell+p-1,
\]
so the desired inequality reduces to
\[
p^\ell-p+1>0.
\]
For the transition $a=\ell-1$ to $a=\ell$, one needs
\[
1+p^\ell>1+\frac{p^{\ell+1}}{1+p^\ell},
\]
which is equivalent to $p^\ell>p-1$.

For $a\ge\ell$, writing $C=p^\ell-p+1$ gives
\[
R_\ell(p,a)=\frac{p^{a+2}-C}{p^{a+1}-C}.
\]
Hence
\[
R_\ell(p,a)-p
=\frac{(p-1)C}{p^{a+1}-C}>0,
\]
and the limit is $p$.
\end{proof}

Since $\log_p$ is strictly increasing, $G_\ell(p,a)=\log_p R_\ell(p,a)$ is also strictly decreasing for integers $a\ge\ell-1$.

For a function $g:S\to\R$, we use the notation
\[
\operatorname*{arg\,max}_{a\in S}g(a)
:=\{a\in S:g(a)\ge g(b)\text{ for every }b\in S\}.
\]
This is the set of arguments at which $g$ attains its maximum.  It is empty if no maximum is attained and may contain several elements when different arguments give the same maximum.  Throughout the following optimization problems, the exponent variables range over integers.

\begin{theorem}[Exact local factorization of the global maximum]\label{thm:localargmax}
Fix $\ell\ge2$ and $\varepsilon>0$.  For each prime $p$, let
\[
\mathcal M_p(\varepsilon):=
\operatorname*{arg\,max}_{a\ge0}f_{p,\varepsilon}(a).
\]
Then:
\begin{enumerate}[label=\textup{(\roman*)}]
\item $\mathcal M_p(\varepsilon)$ is nonempty and finite for every prime $p$;
\item $\mathcal M_p(\varepsilon)=\{0\}$ for all sufficiently large primes $p$;
\item an integer $N$ maximizes $F_{\ell,\varepsilon}$ if and only if
\[
\nu_p(N)\in\mathcal M_p(\varepsilon)
\qquad\text{for every prime }p;
\]
\item if $a\ge\ell$ belongs to $\mathcal M_p(\varepsilon)$, then
\[
R_\ell(p,a-1)\ge p^{1+\varepsilon}\ge R_\ell(p,a).
\]
\end{enumerate}
\end{theorem}

\begin{proof}
For $1\le a<\ell$,
\[
f_{p,\varepsilon}(a)=p^{-a(1+\varepsilon)}<1=f_{p,\varepsilon}(0).
\]
For $a\ge\ell$, Lemma~\ref{lem:local} yields
\[
f_{p,\varepsilon}(a)
\le \frac{p^{a+1}}{(p-1)p^{a(1+\varepsilon)}}
=\frac{p}{p-1}p^{-a\varepsilon},
\]
which tends to $0$ as $a\to\infty$.  Choose an integer $B\ge\ell$ such that
\[
\frac{p}{p-1}p^{-B\varepsilon}<1.
\]
For every $a\ge B$, the preceding bound gives
\[
f_{p,\varepsilon}(a)
\le\frac{p}{p-1}p^{-a\varepsilon}
\le\frac{p}{p-1}p^{-B\varepsilon}
<1=f_{p,\varepsilon}(0).
\]
Thus no exponent $a\ge B$ can maximize $f_{p,\varepsilon}$.  On the finite nonempty set $\{0,1,\ldots,B-1\}$, the function attains a largest value at some exponent $a_*$.  This value is at least $f_{p,\varepsilon}(0)=1$, so it also exceeds every value with $a\ge B$.  Hence $a_*$ is a maximizer over all nonnegative integers, and
\[
\varnothing\ne\mathcal M_p(\varepsilon)
\subseteq\{0,1,\ldots,B-1\}.
\]
This proves both attainment of the local maximum and finiteness of the maximizing set, as asserted in (i).

Also, for $a\ge\ell$,
\[
f_{p,\varepsilon}(a)
\le \frac{p}{p-1}p^{-\ell\varepsilon}.
\]
The right side tends to $0$ as $p\to\infty$.  Hence for all sufficiently large $p$, every positive exponent gives a value below $1=f_{p,\varepsilon}(0)$, proving (ii).

Now let $N=\prod p^{a_p}$.  Since
\[
F_{\ell,\varepsilon}(N)=\prod_p f_{p,\varepsilon}(a_p),
\]
if some $a_p$ is not a local maximizer, replacing only that exponent by an element of $\mathcal M_p(\varepsilon)$ strictly increases the global product.  Conversely, choosing a local maximizer at every prime gives a product at least as large as that of any integer, term by term.  By (ii), only finitely many chosen exponents can be nonzero, so this product corresponds to an integer.  This proves (iii).

Finally, if $a\ge\ell$ is a local maximizer, then
\[
f_{p,\varepsilon}(a)\ge f_{p,\varepsilon}(a-1),
\qquad
f_{p,\varepsilon}(a)\ge f_{p,\varepsilon}(a+1).
\]
Using
\[
\frac{f_{p,\varepsilon}(a)}{f_{p,\varepsilon}(a-1)}
=R_\ell(p,a-1)p^{-(1+\varepsilon)}
\]
and
\[
\frac{f_{p,\varepsilon}(a+1)}{f_{p,\varepsilon}(a)}
=R_\ell(p,a)p^{-(1+\varepsilon)},
\]
we obtain the sandwich inequalities in (iv).
\end{proof}

Parts~\textup{(i)--(iii)} imply that a global maximizer exists for every $\varepsilon>0$: choose one local maximizing exponent at each prime, noting that only finitely many are nonzero.  Every such maximizer is $\ell$-full, because exponents $1,\ldots,\ell-1$ give local values strictly below $f_{p,\varepsilon}(0)=1$.

\section{$\ell$-critical parameters and explicit exponent selection}

Fix $\ell\ge2$, a prime $p$, and $\varepsilon>0$, and put
\[
C_\ell(p):=p^\ell-p+1.
\]
For $a\ge\ell$,
\[
R_\ell(p,a)
=\frac{p^{a+2}-C_\ell(p)}{p^{a+1}-C_\ell(p)}.
\]
Solving $R_\ell(p,a)=p^{1+\varepsilon}$ for a real variable $a$ leads to
\begin{equation}\label{eq:T}
T_\ell(p,\varepsilon)
:=
\log_p\!\left(
\frac{(p^\ell-p+1)(p^{1+\varepsilon}-1)}
{p^2(p^\varepsilon-1)}
\right).
\end{equation}

\begin{theorem}[Positive-branch maximizer]\label{thm:positivebranch}
Let
\[
\mathcal M_p^+(\varepsilon)
:=\operatorname*{arg\,max}_{a\ge\ell}f_{p,\varepsilon}(a),
\]
and write $T=T_\ell(p,\varepsilon)$.  Then
\[
\mathcal M_p^+(\varepsilon)=
\begin{cases}
\{\ell\},&T<\ell,\\[1mm]
\{m,m+1\},&T=m\in\mathbb Z,\ m\ge\ell,\\[1mm]
\{\lfloor T\rfloor+1\},&T>\ell\text{ and }T\notin\mathbb Z.
\end{cases}
\]
\end{theorem}

\begin{proof}
Write $f(a)=f_{p,\varepsilon}(a)$, $C=p^\ell-p+1$, and $q=p^{1+\varepsilon}$.  Then $q>p$ and $p^{a+1}-C>0$ for every $a\ge\ell$.  For such $a$,
\[
\frac{f(a+1)}{f(a)}
=\frac{R_\ell(p,a)}q.
\]
Using the formula for $R_\ell(p,a)$ and the positivity of the denominators, we obtain
\begin{align*}
\frac{f(a+1)}{f(a)}>1
&\iff p^{a+2}-C>q\bigl(p^{a+1}-C\bigr)\\
&\iff (q-p)p^{a+1}<C(q-1)\\
&\iff p^a<\frac{C(q-1)}{p(q-p)}\\
&\iff p^a<\frac{(p^\ell-p+1)(p^{1+\varepsilon}-1)}
{p^2(p^\varepsilon-1)}=p^T\\
&\iff a<T.
\end{align*}
The same rearrangement with equality or with the reversed inequality shows that
\[
\frac{f(a+1)}{f(a)}
\begin{cases}
>1,&a<T,\\
=1,&a=T,\\
<1,&a>T.
\end{cases}
\]

If $T<\ell$, then every integer $a\ge\ell$ satisfies $a>T$.  Hence $f(a+1)<f(a)$ throughout the positive-exponent branch.  Its unique maximum is therefore attained at the first exponent $a=\ell$, giving $\mathcal M_p^+(\varepsilon)=\{\ell\}$.

If $T=m\in\mathbb Z$ with $m\ge\ell$, then $f(a+1)>f(a)$ for $\ell\le a<m$, $f(m+1)=f(m)$, and $f(a+1)<f(a)$ for $a>m$.  Thus every exponent $\ell\le a<m$ gives a value smaller than $f(m)$, and every exponent above $m+1$ gives a value smaller than $f(m+1)$.  The maximum is attained exactly at $m$ and $m+1$.  This also covers $m=\ell$, when there are no exponents in the increasing part before the tie.

Finally, suppose $T>\ell$ and $T\notin\mathbb Z$, and put $k=\lfloor T\rfloor$.  Then $k\ge\ell$ and $k<T<k+1$.  For every integer $\ell\le a\le k$, we have $f(a+1)>f(a)$, whereas for every integer $a\ge k+1$, we have $f(a+1)<f(a)$.  The sequence therefore increases strictly up to exponent $k+1$ and decreases strictly thereafter.  Its unique maximizing exponent is $k+1=\lfloor T\rfloor+1$.
\end{proof}

The positive-exponent sequence $(f_{p,\varepsilon}(a))_{a\ge\ell}$ is therefore \emph{unimodal}: it increases strictly up to its maximum, may take that maximum at two adjacent exponents, and then decreases strictly.  The increasing part may be empty when the maximum occurs at the initial exponent $\ell$.  Thus the maximum occurs either once or at precisely two consecutive exponents.

\begin{theorem}[Complete global exponent rule]\label{thm:globalrule}
For a prime $p$, let $\mathcal M_p^+(\varepsilon)$ be as in Theorem~\ref{thm:positivebranch}, and let
\[
V_p(\varepsilon):=
\max_{a\ge\ell}f_{p,\varepsilon}(a).
\]
Then the full local argmax is
\[
\mathcal M_p(\varepsilon)=
\begin{cases}
\{0\},&V_p(\varepsilon)<1,\\[1mm]
\mathcal M_p^+(\varepsilon),&V_p(\varepsilon)>1,\\[1mm]
\{0\}\cup\mathcal M_p^+(\varepsilon),&V_p(\varepsilon)=1.
\end{cases}
\]
\end{theorem}

\begin{proof}
For the fixed prime $p$, abbreviate $f(a)=f_{p,\varepsilon}(a)$ and $V=V_p(\varepsilon)$.  At exponent $0$ we have $f(0)=1$.  If $1\le a<\ell$, then $\sigma_{(\ell)}(p^a)=1$, so
\[
f(a)=p^{-a(1+\varepsilon)}<1.
\]
Consequently, none of the exponents $1,\ldots,\ell-1$ can maximize $f$ over all $a\ge0$.

By Theorem~\ref{thm:positivebranch}, the maximum $V$ over $a\ge\ell$ is attained, and
\[
f(a)\le V\quad(a\ge\ell),
\qquad
f(a)=V\iff a\in\mathcal M_p^+(\varepsilon).
\]
We compare this value with $f(0)=1$.

If $V<1$, then every $a\ge\ell$ satisfies $f(a)\le V<1$, and the exponents $1,\ldots,\ell-1$ also give values below $1$.  Thus $a=0$ is the unique maximizer, and $\mathcal M_p(\varepsilon)=\{0\}$.

If $V>1$, then exponent $0$ and all exponents below $\ell$ give values strictly smaller than $V$.  Among the exponents $a\ge\ell$, precisely those in $\mathcal M_p^+(\varepsilon)$ attain $V$.  Therefore $\mathcal M_p(\varepsilon)=\mathcal M_p^+(\varepsilon)$.

If $V=1$, then exponent $0$ and every exponent in $\mathcal M_p^+(\varepsilon)$ attain the value $1$.  Every exponent $1\le a<\ell$ gives a smaller value, and any $a\ge\ell$ outside $\mathcal M_p^+(\varepsilon)$ satisfies $f(a)<V=1$.  Hence the complete maximizing set is $\{0\}\cup\mathcal M_p^+(\varepsilon)$.
\end{proof}

Combining this rule with Theorem~\ref{thm:localargmax}, we obtain the global parametrization: all $\ell$-colossally abundant numbers for the parameter $\varepsilon$ are exactly the integers
\[
N=\prod_p p^{a_p}
\qquad\text{with}\qquad
a_p\in\mathcal M_p(\varepsilon)
\]
for every prime $p$.  In particular, the maximizer is unique if and only if every $\mathcal M_p(\varepsilon)$ is a singleton.

\begin{remark}[Two distinct kinds of ties]\label{rem:ties}
In the classical case $\ell=1$, the possible exponents $0,1,2,\ldots$ form a single consecutive sequence.  The notation $0\leftrightarrow1$ denotes a change in the maximizing exponent of a prime $p$ between $0$ and $1$: exponent $0$ means that $p$ is absent from the integer, whereas exponent $1$ means that $p$ divides it exactly once.  At this transition both exponents give the same maximal local value,
\[
f_{p,\varepsilon}(0)=1=f_{p,\varepsilon}(1)
=\frac{1+p}{p^{1+\varepsilon}},
\]
so the transition occurs at $\varepsilon=\log_p(1+1/p)$.  The classical tied choices are described in \cite{CaveneyEtAl2012}.

For $\ell\ge2$, the exponents $1,\ldots,\ell-1$ are always dominated by exponent $0$, and a prime may enter by jumping directly from exponent $0$ to an exponent at least $\ell$.  We therefore distinguish two mechanisms for nonuniqueness.
\begin{enumerate}[label=\textup{(\roman*)}]
\item \emph{Adjacent positive-exponent ties:} for $a\ge\ell$,
\[
\varepsilon=\log_p R_\ell(p,a)-1
\]
forces $f_{p,\varepsilon}(a)=f_{p,\varepsilon}(a+1)$.  This affects the global maximizer only when their common value is at least $1$.
\item \emph{Entry ties with exponent $0$:} for $a\ge\ell$,
\[
\varepsilon=\frac{\log_p\sigma_{(\ell)}(p^a)}{a}-1
\]
forces $f_{p,\varepsilon}(a)=1=f_{p,\varepsilon}(0)$.  It contributes to nonuniqueness of the global maximum only when $a$ is a positive-branch maximizing exponent at that same parameter.
\end{enumerate}
Thus adjacent-ratio equalities alone do not describe all changes of the global maximizer.
\end{remark}

\begin{definition}[Candidate and global critical parameters]\label{def:critical}
For $a\ge\ell$ define the adjacent-transition value
\[
\varepsilon^{\mathrm{adj}}_{\ell,p,a}
:=\log_p R_\ell(p,a)-1>0,
\]
and define the entry value
\[
\varepsilon^{\mathrm{ent}}_{\ell,p,a}
:=\frac{\log \sigma_{(\ell)}(p^a)}{a\log p}-1.
\]
Only positive entry values are relevant.  Put
\[
\mathfrak E_\ell
:=
\{\varepsilon^{\mathrm{adj}}_{\ell,p,a}:p\ \text{prime},\ a\ge\ell\}
\cup
\{\varepsilon^{\mathrm{ent}}_{\ell,p,a}>0:p\ \text{prime},\ a\ge\ell\}.
\]
We call $\varepsilon>0$ a \emph{global $\ell$-critical parameter} if $F_{\ell,\varepsilon}$ has more than one global maximizer.
\end{definition}

\begin{theorem}[Uniqueness outside the candidate critical set]\label{thm:criticalset}
The set $\mathfrak E_\ell$ is countable.  If $\varepsilon\notin\mathfrak E_\ell$, then $F_{\ell,\varepsilon}$ has a unique global maximizer.  For every fixed $\varepsilon>0$, the set of global maximizers is finite.
\end{theorem}

\begin{proof}
Countability is immediate because $\mathfrak E_\ell$ is indexed by pairs $(p,a)$ with $p$ prime and $a\ge\ell$ integer.

Fix $\varepsilon\notin\mathfrak E_\ell$.  For every prime $p$, no equality
\[
R_\ell(p,a)=p^{1+\varepsilon}
\]
can occur at an integer $a\ge\ell$.  Hence Theorem~\ref{thm:positivebranch} gives a unique positive-branch maximizing exponent $a_p^+$.  If its value were equal to $1$, then
\[
\frac{\sigma_{(\ell)}(p^{a_p^+})}{p^{a_p^+(1+\varepsilon)}}=1,
\]
which rearranges to
\[
\varepsilon
=\frac{\log\sigma_{(\ell)}(p^{a_p^+})}{a_p^+\log p}-1
=\varepsilon^{\mathrm{ent}}_{\ell,p,a_p^+},
\]
contrary to $\varepsilon\notin\mathfrak E_\ell$.  Therefore the comparison in Theorem~\ref{thm:globalrule} is strict for every prime, and each local argmax $\mathcal M_p(\varepsilon)$ is a singleton.  The global maximizer is consequently unique.

For arbitrary fixed $\varepsilon>0$, Theorem~\ref{thm:localargmax} shows that each local argmax is finite and that $\mathcal M_p(\varepsilon)=\{0\}$ for all sufficiently large $p$.  Hence only finitely many primes contribute more than one local choice, and the Cartesian product of the local choice sets is finite.  By Theorem~\ref{thm:globalrule}, this Cartesian product parametrizes all global maximizers.
\end{proof}

In particular, every global $\ell$-critical parameter belongs to $\mathfrak E_\ell$.

\begin{remark}\label{rem:candidate-not-global}
For every $\ell\ge2$, membership in $\mathfrak E_\ell$ is necessary but not sufficient for $\varepsilon$ to be a global $\ell$-critical parameter.
\end{remark}

\begin{proof}
Necessity follows from Theorem~\ref{thm:criticalset}.  To disprove sufficiency, set
\[
\varepsilon_0
:=\log_2\!\left(\frac{1+3\cdot2^\ell}{2(1+2^\ell)}\right)
=\varepsilon^{\mathrm{adj}}_{\ell,2,\ell}\in\mathfrak E_\ell.
\]
Since $2^\ell\ge4$, we have
\[
2^{\varepsilon_0}\ge\frac{13}{10},
\qquad
\left(\frac{13}{10}\right)^4>2,
\]
and hence $\varepsilon_0>1/4$.  At the prime $2$, Theorem~\ref{thm:positivebranch} gives the positive maximizing exponents $\ell$ and $\ell+1$.  Their common value satisfies
\[
f_{2,\varepsilon_0}(\ell)
=f_{2,\varepsilon_0}(\ell+1)
=(1+2^{-\ell})2^{-\ell\varepsilon_0}
<\frac{5}{4\sqrt2}<1.
\]
Thus $\mathcal M_2(\varepsilon_0)=\{0\}$.  For every prime $p\ge3$ and integer $a\ge\ell$,
\[
f_{p,\varepsilon_0}(a)
\le\frac{p}{p-1}p^{-a\varepsilon_0}
<\frac{\sqrt p}{p-1}
\le\frac{\sqrt3}{2}<1.
\]
The second inequality uses $a\ge2$ and $\varepsilon_0>1/4$; the third follows from the strict decrease of $x\mapsto\sqrt{x}/(x-1)$ for $x>1$.  Exponents $1,\ldots,\ell-1$ also give values below $1$.  Consequently $\mathcal M_p(\varepsilon_0)=\{0\}$ for every prime $p$, and Theorem~\ref{thm:globalrule} gives the unique global maximizer $N=1$.  Therefore $\varepsilon_0$ is not a global $\ell$-critical parameter.
\end{proof}

An adjacent positive-exponent equality contributes to global nonuniqueness only when its common value is at least $1$.  Likewise, an entry equality $f_{p,\varepsilon}(a)=1$ contributes only when $a$ maximizes the positive branch at that parameter.

\begin{proposition}[Geometry of the candidate critical set]\label{prop:criticalgeometry}
Fix $\ell\ge2$.  For every $\delta>0$, the set
\[
\mathfrak E_\ell\cap[\delta,\infty)
\]
is finite.
\end{proposition}

\begin{proof}
We treat the two families in Definition~\ref{def:critical} separately.

For the adjacent-transition values, Lemma~\ref{lem:Rdecreasing} gives, for every $a\ge\ell$,
\[
R_\ell(p,a)\le R_\ell(p,\ell)
=\frac{1+p^\ell+p^{\ell+1}}{1+p^\ell}
<p+1.
\]
Hence
\[
0<\varepsilon^{\mathrm{adj}}_{\ell,p,a}
=\log_p\!\left(\frac{R_\ell(p,a)}p\right)
<\log_p\!\left(1+\frac1p\right),
\]
and the last quantity tends to $0$ as $p\to\infty$.  Thus, for a fixed $\delta>0$, only finitely many primes can satisfy
$\varepsilon^{\mathrm{adj}}_{\ell,p,a}\ge\delta$.  For each such prime, the adjacent-transition values are strictly decreasing in $a$ and tend to $0$ by Lemma~\ref{lem:Rdecreasing}, so only finitely many exponents can contribute values at least $\delta$.

For the entry values, when $a\ge\ell$ we have
\[
1<A_\ell(p^a)<\frac{p}{p-1}.
\]
Therefore
\[
0<\varepsilon^{\mathrm{ent}}_{\ell,p,a}
=\frac{\log_p A_\ell(p^a)}a
<\frac1a\log_p\!\left(\frac{p}{p-1}\right)
\le \frac1\ell\log_p\!\left(\frac{p}{p-1}\right).
\]
The final bound tends to $0$ as $p\to\infty$, so again only finitely many primes can yield an entry value at least $\delta$.  Once $p$ is fixed, the inequality
$\varepsilon^{\mathrm{ent}}_{\ell,p,a}\ge\delta$ implies
\[
a<\frac{\log_p(p/(p-1))}{\delta},
\]
so only finitely many exponents occur.

This proves the asserted finiteness.
\end{proof}

For any fixed prime $p$, the distinct values $\varepsilon^{\mathrm{adj}}_{\ell,p,a}$ decrease to $0$ as $a\to\infty$.  Thus $\mathfrak E_\ell$ is infinite and has $0$ as an accumulation point.  The proposition excludes every positive accumulation point and bounds the number of elements above any positive threshold.  Consequently, after repetitions are removed, its elements can be listed as
\[
\eta_1>\eta_2>\eta_3>\cdots>0,
\qquad \eta_j\longrightarrow0.
\]
Hence $0$ is the only accumulation point of $\mathfrak E_\ell$ in $[0,\infty)$.

\begin{proposition}[Classical specialization of the transition formula]\label{prop:l1T}
Let $\ell=1$ and $\varepsilon>0$.  For each prime $p$,
\[
T_1(p,\varepsilon)
=\log_p\!\left(\frac{p^{1+\varepsilon}-1}{p^2(p^\varepsilon-1)}\right)>-1.
\]
If $T_1(p,\varepsilon)\notin\mathbb Z$, the unique maximizing exponent over all integers $a\ge0$ is
\[
\alpha_p(\varepsilon)=\lfloor T_1(p,\varepsilon)\rfloor+1.
\]
If $T_1(p,\varepsilon)=m\in\mathbb Z$, the two maximizing exponents are $m$ and $m+1$.
\end{proposition}

\begin{proof}
For $\ell=1$, every divisor is admissible, and for every integer $a\ge0$,
\[
\sigma_{(1)}(p^a)=1+p+\cdots+p^a
=\frac{p^{a+1}-1}{p-1}.
\]
Hence
\[
R_1(p,a)=\frac{p^{a+2}-1}{p^{a+1}-1},
\qquad
\frac{f_{p,\varepsilon}(a+1)}{f_{p,\varepsilon}(a)}
=R_1(p,a)p^{-(1+\varepsilon)}.
\]
Write $q=p^{1+\varepsilon}>p$.  Since $p^{a+1}-1>0$, the transition equation is equivalent to
\begin{align*}
R_1(p,a)=q
&\iff p^{a+2}-1=q\bigl(p^{a+1}-1\bigr)\\
&\iff (q-p)p^{a+1}=q-1\\
&\iff p^a=\frac{q-1}{p(q-p)}
=\frac{p^{1+\varepsilon}-1}{p^2(p^\varepsilon-1)}.
\end{align*}
Taking the logarithm to base $p$ gives the stated expression for $T_1=T_1(p,\varepsilon)$.  Moreover,
\[
p^{T_1}
=\frac1{p^2}\left(p+\frac{p-1}{p^\varepsilon-1}\right)
>\frac1p,
\]
so $T_1>-1$.

To determine the maximizing exponents, the same algebra gives
\begin{align*}
\frac{f_{p,\varepsilon}(a+1)}{f_{p,\varepsilon}(a)}>1
&\iff p^{a+2}-1>q\bigl(p^{a+1}-1\bigr)\\
&\iff (q-p)p^{a+1}<q-1\\
&\iff p^a<p^{T_1}\\
&\iff a<T_1.
\end{align*}
Equality holds exactly when $a=T_1$, and the ratio is less than $1$ exactly when $a>T_1$.

Suppose first that $T_1\notin\mathbb Z$, and put $k=\lfloor T_1\rfloor+1$.  Because $T_1>-1$, we have $k\ge0$.  If $k=0$, then $-1<T_1<0$, so every integer $a\ge0$ satisfies $a>T_1$; the sequence is strictly decreasing, and its unique maximum is at $a=0$.  If $k\ge1$, then $k-1<T_1<k$.  The adjacent ratio is greater than $1$ for $0\le a\le k-1$ and less than $1$ for $a\ge k$.  Thus the sequence increases strictly up to $a=k$ and decreases strictly thereafter.  In both cases the unique maximizing exponent is $k=\lfloor T_1\rfloor+1$.

Now suppose $T_1=m\in\mathbb Z$.  The inequality $T_1>-1$ gives $m\ge0$.  The adjacent ratio is greater than $1$ for $0\le a<m$, equals $1$ at $a=m$, and is less than $1$ for $a>m$.  Therefore
\[
f_{p,\varepsilon}(m)=f_{p,\varepsilon}(m+1),
\]
and this common value is strictly larger than the value at every other nonnegative exponent.  The maximizing set is exactly $\{m,m+1\}$, including the boundary case $m=0$, where the sequence begins with the tie at exponents $0$ and $1$ and then decreases strictly.
\end{proof}

Since $T_1>-1$, an integral transition value $T_1=m$ necessarily has $m\ge0$.  The identity
\[
T_1(p,\varepsilon)
=\log_p\!\left(\frac{p^{1+\varepsilon}-1}{p^\varepsilon-1}\right)-2
\]
gives the equivalent formula, away from ties,
\[
\alpha_p(\varepsilon)
=\left\lfloor
\log_p\!\left(\frac{p^{1+\varepsilon}-1}{p^\varepsilon-1}\right)
\right\rfloor-1.
\]
Away from ties, this agrees with the exponent formula of Alaoglu and Erd\H{o}s \cite[Theorem~10]{AlaogluErdos1944}.  At a tie, our nonstrict definition retains both adjacent exponents, as in the classical treatment of \cite{CaveneyEtAl2012}.

\begin{remark}[Simultaneous transitions]\label{rem:simultaneous}
For classical colossally abundant numbers, coincidences between prime transition parameters and their connection with transcendence theory are discussed in \cite{CaveneyEtAl2012}; see also Lang \cite{Lang1966} for general background.  Theorem~\ref{thm:globalrule} describes all simultaneous ties in the restricted setting through the local maximizing sets, without requiring a bound on the number of primes involved.
\end{remark}

\begin{theorem}\label{thm:infCA}
For every fixed $\ell\ge2$, there are infinitely many $\ell$-colossally abundant numbers.
\end{theorem}

\begin{proof}
Fix a prime $p$.  For every integer $a\ge\ell$, define
\[
\varepsilon_a:=\log_p R_\ell(p,a)-1.
\]
By Lemma~\ref{lem:Rdecreasing},
\[
\varepsilon_a>0,
\qquad
\varepsilon_a\downarrow0.
\]
At $\varepsilon=\varepsilon_a$ we have
\[
f_{p,\varepsilon_a}(a)=f_{p,\varepsilon_a}(a+1),
\]
and Theorem~\ref{thm:positivebranch} shows that these two exponents maximize the positive branch.

It remains to prove that their common value exceeds $1$ for all sufficiently large $a$.  Put $C=p^\ell-p+1$ and
\[
u_a:=\frac{C}{p^{a+1}}=\frac Cp\,p^{-a}.
\]
The coefficient $C/p$ depends only on the fixed parameters $p$ and $\ell$, and is independent of $a$.  Thus $u_a=O_{p,\ell}(p^{-a})$ and $u_a\to0$.  Also, $0<C<p^\ell$ and $a\ge\ell$ imply $0<u_a<1/p<1$.  We have
\[
\frac{R_\ell(p,a)}p
=\frac{1-u_a/p}{1-u_a}
=1+\frac{(1-1/p)u_a}{1-u_a}.
\]
Using $\log(1+t)\le t$ for $t\ge0$, we obtain
\[
0<\varepsilon_a
=\log_p\!\left(1+\frac{(1-1/p)u_a}{1-u_a}\right)
\le
\frac{(1-1/p)u_a}{(1-u_a)\log p}.
\]
For all sufficiently large $a$, we have $u_a\le1/2$.  The preceding estimate and the identity $u_a=(C/p)p^{-a}$ then give
\[
0<a\varepsilon_a
\le\frac{2(1-1/p)C}{p\log p}\,a p^{-a}
\longrightarrow0.
\]
Here $a p^{-a}\to0$ because the ratio of consecutive terms is $(a+1)/(ap)$, which tends to $1/p<1$.
Now
\[
f_{p,\varepsilon_a}(a)
=A_\ell(p^a)p^{-a\varepsilon_a}.
\]
By Lemma~\ref{lem:primepowermonotone},
\[
A_\ell(p^a)\longrightarrow\frac{p}{p-1},
\]
while $p^{-a\varepsilon_a}\to1$.  Hence
\[
f_{p,\varepsilon_a}(a)
\longrightarrow\frac{p}{p-1}>1.
\]
Thus for all sufficiently large $a$, both $a$ and $a+1$ belong to the full local argmax $\mathcal M_p(\varepsilon_a)$.

For each such $a$, choose a global maximizer $N_a$ with $\nu_p(N_a)=a$; Theorem~\ref{thm:globalrule} guarantees that such a choice is possible.  The integers $N_a$ are distinct because their $p$-adic valuations are distinct.  Hence there are infinitely many $\ell$-colossally abundant numbers.
\end{proof}

\section{Explicit examples}

The preceding theorems describe every $\ell$-CA maximizer for a fixed $\varepsilon$ through finitely many local comparisons: determine the positive-branch maximizing exponents from $T_\ell(p,\varepsilon)$, compare their common value with $1$, and use the large-prime bound to exclude all remaining primes.  The following examples illustrate this description; the decisive comparisons are verified by exact integer inequalities.  In these examples, a prime is called active when its unique local maximizing exponent is positive.

\begin{example}\label{ex:l2}
Let $\ell=2$.

\smallskip
\noindent\textup{(a) $\varepsilon=1/5$.}
For $p=2$, formula \eqref{eq:T} gives
\[
T_2(2,1/5)\approx2.710121594,
\]
so the positive-branch maximizing exponent is $3$, and
\[
f_{2,1/5}(3)=\frac{13}{2^{18/5}}\approx1.072100178>1.
\]
The choice of exponent $3$ and the comparison with $1$ can be checked exactly:
\[
13^5>5^5\,2^6,
\qquad
29^5<13^5\,2^6,
\qquad
13^5>2^{18}.
\]
These inequalities give $f_{2,1/5}(3)>f_{2,1/5}(2)$, $f_{2,1/5}(3)>f_{2,1/5}(4)$, and $f_{2,1/5}(3)>1$, respectively.  Unimodality therefore shows that $2$ is active with the unique exponent $3$.

For every prime $p\ge3$ and every $a\ge2$,
\[
f_{p,1/5}(a)
\le\frac{p}{p-1}p^{-2/5}.
\]
Writing the bound as $p^{3/5}/(p-1)$, its logarithmic derivative is $\frac{3}{5p}-\frac1{p-1}<0$ for $p>1$; at $p=3$ its fifth power is $27/32<1$.  Hence all odd primes are inactive.  Therefore the unique maximizer is
\[
N_{1/5}=2^3=8.
\]

\smallskip
\noindent\textup{(b) $\varepsilon=1/100$.}
The local data are
\[
\begin{array}{c|c|c|c}
 p&T_2(p,1/100)&\text{positive maximizer }a_p^+&f_{p,1/100}(a_p^+)\\
\hline
2&6.772513337&7&1.882948543\\
3&4.518328553&5&1.406192786\\
5&3.326158757&4&1.164187424\\
7&2.884768407&3&1.080800381\\
11&2.486028356&3&1.015887878\\
13&2.374078289&3&0.997584977\\
17&2.221953921&3&0.972732876
\end{array}
\]
For an exact verification, put $S_p(a)=\sigma_{(2)}(p^a)$.  Every row satisfies the integer inequalities
\[
S_p(a_p^+)^{100}>p^{101}S_p(a_p^+-1)^{100},
\qquad
p^{101}S_p(a_p^+)^{100}>S_p(a_p^++1)^{100}.
\]
Thus the local value at the displayed exponent strictly exceeds both neighboring local values, and unimodality proves that it is the unique positive maximizer.  Moreover,
\[
S_p(a_p^+)^{100}
\begin{cases}
>p^{101a_p^+},&p\in\{2,3,5,7,11\},\\
<p^{101a_p^+},&p\in\{13,17\}.
\end{cases}
\]
Hence $2,3,5,7,11$ are active, while $13$ and $17$ are not.  For $p\ge19$ and $a\ge2$,
\[
f_{p,1/100}(a)
\le\frac{p}{p-1}p^{-1/50}.
\]
Writing the bound as $p^{49/50}/(p-1)$, its logarithmic derivative is $\frac{49}{50p}-\frac1{p-1}<0$; at $p=19$ the bound is below $1$ because the exact integer comparison $19^{49}<18^{50}$ holds.  Hence every prime $p\ge19$ is inactive.  Consequently the unique maximizer is
\[
N_{1/100}=2^7\,3^5\,5^4\,7^3\,11^3.
\]
\end{example}

\section{Conclusion}

The divisor restriction studied in this paper preserves enough multiplicative structure for a substantial part of the classical abundance theory to survive, but it changes the local arithmetic in a decisive way.  At $\ell=1$ we recover the ordinary sum of divisors and hence the usual abundant, perfect, superabundant, and colossally abundant theories.  For $\ell\ge2$, the gap in the allowed prime exponents produces a strict hierarchy of abundance classes and forces every nontrivial $\ell$-superabundant and $\ell$-colossally abundant number to be $\ell$-full.

The superabundant results extend the classical consecutive-prime and nonincreasing-exponent structure, while the colossally abundant results extend the Ramanujan--Alaoglu--Erd\H{o}s prime-by-prime optimization principle.  The explicit transition formula
\[
T_\ell(p,\varepsilon)
=\log_p\!\left(
\frac{(p^\ell-p+1)(p^{1+\varepsilon}-1)}
{p^2(p^\varepsilon-1)}
\right)
\]
provides a closed description of the positive local maximizer.  The additional comparison with exponent $0$ is essential and yields the complete global rule, including every possible tie.  This characterization also gives a direct proof that infinitely many $\ell$-colossally abundant numbers exist for each fixed $\ell\ge2$.

Together, these results extend the structural and extremal theory of classical abundance to the family of $\ell$-restricted divisor sums.

\end{document}